\documentclass[12pt]{amsart}
\usepackage{amsmath,amscd,latexsym,verbatim,amssymb}
\usepackage{times}       
 
\newcommand{\cf}{{\it cf.} }
\newcommand{\ie}{{\it i.e.} }
\newcommand{\eg}{{\it e.g.} }

\newcommand{\Q}{\mathbb{Q}}

\newcommand{\C}{\mathbb{C}}

\newcommand{\sK}{{\mathcal{K}}}

\newcommand{\inj}{\hookrightarrow}

 \newcommand{\e}{\frac{1}{p^\infty}}

\newcommand{\Spec}{\operatorname{Spec}}

 \newcounter{spec}

\swapnumbers

\newtheorem{thm}{Theorem}[subsection]
\newtheorem{lemma}[thm]{Lemma}
\newtheorem{prop}[thm]{Proposition}

\theoremstyle{definition}
\newtheorem{defn}[thm]{Definition}

\numberwithin{equation}{section}

 at10pt
 at12pt

\begin{document}

\title[Singularities in mixed characteristic. The perfectoid approach.]{Singularities in mixed characteristic. The perfectoid approach.}
 
\author{Y. Andr\'e}
 
  \address{Sorbonne-Universit\'e, Institut de Math\'ematiques de Jussieu, 4 place Jussieu, 75005 Paris France.}
\email{yves.andre@imj-prg.fr}
  \begin{abstract}  The homological conjectures, which date back to Peskine, Szpiro and Hochster in the late 60Õs, make fundamental predictions about syzygies and intersection problems in commutative algebra. They were settled long ago in the presence of a base field and led to tight closure theory, a powerful tool to investigate singularities in characteristic $p$. 

Recently, perfectoid techniques coming from $p$-adic Hodge theory have allowed us to get rid of any base field; this solves the direct summand conjecture and establishes the existence and weak functoriality of big Cohen-Macaulay algebras, which solve in turn the homological conjectures in general. 
 This also opens the way to the study of singularities in mixed characteristic. 
 
 We sketch a broad outline of this story, taking lastly a glimpse at ongoing work by L. Ma and K. Schwede, which shows how such a study could build a bridge between singularity theory in char. $p$ and in char. $0$. 
  \end{abstract}
\maketitle

 \begin{sloppypar}   
  
 \bigskip
 
 \bigskip

 \section{Commutative Algebra is not (only) Chapter 0 of Algebraic Geometry.}
 
  \subsection{} Commutative Algebra took its roots in Algebraic Number Theory, Algebraic Geometry and Invariant Theory at the end of the XIXth century, in the wake of Dedekind, Kronecker and Gordan. It soon provided a unified language for the three theories as well as a consolidation of their foundations. 
  
  It may be argued that Commutative Algebra started with Hilbert's three fundamental theorems (which still bear their original names: {\it basis theorem, syzygy theorem, Nullstellensatz}), and was put on firm ground in the 30's with Krull's four cardinal concepts ({\it locality, dimension, completion and regularity}). 
  
 A main source of examples and inspiration is given by affine rings (rings of functions on affine algebraic varieties over some base field). {\it Regular local rings} (those having the property that the minimal number of generators of their maximal ideal is their Krull dimension) then generalize rings of functions around a {\it nonsingular} point of an algebraic variety. Following this viewpoint, general {\it local rings} appear as algebraic counterparts of {\it singularities} (vs. nonsingular points). 
 
   The theory was so successful that its relation to Algebraic Geometry changed: it was the latter to be founded on the former, in Zariski's and Grothendieck's time. Throughout the canonical text
 EGA, Commutative Algebra forms an extended, crawling ``Chapter $0$". 
  
\subsection{}  However, Commutative Algebra is more than the solfeggio of the music of Algebraic Geometry. 
 In fact, already in the late 50's, a turn happened in the story of Commutative Algebra: {\it the homological turn}, which  was impulsed by Auslander, Buchsbaum and Serre, and originated in Hilbert's syzygy's theorem (and also from intersection theory and other considerations). From general theory of (Noetherian) commutative rings and their ideals, Commutative Algebra became the homological theory of (Noetherian) modules.
 
  In this new viewpoint, regular local rings are those for which any finite module has a free resolution of finite length.  
 
  By the way, let us recall that syzygies generalize the notion of presentation of a module. If a finite module $M$ over a local ring $S$ admits a finite free resolution $0\to S^{b_s}\to S^{b_{s-1}}\to \ldots \to S^{b_0} \to M\to 0$, which we assume to be minimal, the $b_i$'s are called the {\it Betti numbers} of $M$ (so that $b_0$ is the minimal number of generators, $b_1$ the minimal number of relations between them,...), and Im$\,S^{b_i}$ is called the {\it $i$-th syzygy} of $M$. Describing syzygies and estimating Betti numbers is a traditional and typical task of Commutative Algebra, not necessarily related to any specific problem from Algebraic Geometry.

   \subsection{} An important leimotiv in Commutative Algebra is that singularities are always ``finite above non-singularities". One version is Noether's theorem: every affine algebra is finite over a polynomial algebra. Another version, in the complete local case, is Cohen's theorem: every complete local domain is finite over a regular complete local ring (in fact a ring of formal power series with coefficients in a field or a discrete valuation ring). 
   
 The study of singularities thus focusses on {\it finite extensions of a regular (local) ring $S$}. $S$-modules need not have finite free resolutions; in fact, their homological study is complicated and is the subject of the {\it homological conjectures} which we will touch upon. But as we will see, it is nevertheless useful to represent the ring $S$ as a finite extension of a regular ring $R$.

 \subsection{}
 We shall report on 
 recent progresses which came from an encounter between two domains:
 $\;$

\bigskip

\noindent 
{$\,$ Commutative Algebra} $------- $ {$\,$ $p$-adic Hodge Theory}

$\;\;\;$ (Hilbert, Krull, ...) $\;\;\;\;\;\;\;\;\;\;\;\;\;\;\;\;\;\;\;\;\;\;\;\;\;\;\;\;$ (Tate, Fontaine, Faltings...)

$\;\;\;$ Noetherian world  $\;\;\;\;\;\;\;\;\;\;\;\;\;\;\;\;\;\;\;\;\;\;\;\;\;\;\;\;\;\;$ non-Noetherian world 

(finite-dimensional rings, $\;\;\;\;\;\;\;\;\;\;\;\;\;\;\;\;\;\;\;\;\;\;$ (non-archimedean  

$\;\;\;\;\;\;\;\;\;$  finite modules)  $\;\;\;\;\;\;\;\;\;\;\;\;\;\;\;\;\;\;\;\;\;\;\;\;\;\;\;\;\;$ Banach algebras...) 

 \bigskip 

 {$\,$ Homological conjectures} $\;\;\;\;\stackrel{new}{\leftarrow} ---- $ {$\,$ perfectoid theory}

(Peskine-Szpiro, Hochster...) $\;\;\;\;\;\;\;\;\;\;\;\;\;\;\;\;\;\;\;\;$ (Scholze, ...)

\bigskip The second domain, $p$-adic Hodge Theory, which started with Tate's work in the late 60s, studies $p$-adic representations of $p$-adic Galois groups, especially those which come from the cohomology of algebraic varieties. Perfectoid theory is a new development ($\sim 2010$) which has already profoundly transformed $p$-adic Hodge Theory. 

In the ``Chapter $0$" viewpoint, Commutative Algebra stands as the provider of basic concepts, tools and theorems for Algebraic or Arithmetic Geometry and its advanced developments, such as $p$-adic Hodge Theory. 

In the present story, we will see $p$-adic Hodge Theory contribute backward to some basic problems of Commutative Algebra.

\section{The direct summand theorem.}

\subsection{} Let us tell the first instance of this encounter.

 \smallskip Let $R$ be a 
 Noetherian commutative ring, and $R\subset S$ be a finite extension, \ie $S$ is a faithful $R$-algebra, finitely generated as $R$-module (as in Noether's or Cohen's theorem). One thus has an exact sequence of 
  finite $R$-modules:
$$(*)\;\;\;0\to R\to S \to S/R\to 0.$$
 
\smallskip\noindent {
{\bf Question 1}: {\it does this exact sequence 
{splits}?}}
 
\smallskip\noindent Equivalently: is $R$ a direct summand in $S$?
Is there an $R$-linear map $\lambda: S\to R$ with $\lambda(1)= 1$?

\medskip This basic question arose in the late 60's (Hochster \cite{H1}, Raynaud-Gruson \cite{RG}), in the framework of the homological conjectures and of descent theory (independently). An algebraic geometer may not find easy to confer it any geometric meaning; this illustrates what we wrote about ``Chapter $0$". But here is a related question with a more immediate algebro-geometric content:

 \smallskip\noindent  
{\bf Question 2}: {\it is there an $S$-algebra $T$ which is faithfully flat over $R$?}
 
 \smallskip\noindent Equivalently: is $\Spec S \to \Spec R$ a covering for the fpqc topology? 
 
 (It is not required that $T$ is finitely generated, nor even noetherian).

\medskip Actually, a positive answer to Question 2 implies a positive answer to Question 1. Indeed: 

 $a)$ any faithfully flat map $R\to T$ is {\it pure}, \ie universally injective.

$b)$ if a composition $R\to S\to T$ is pure, so is $R\to S$,

$c)$ a finite extension $R\subset S$ is pure if and only if it splits (in the $R$-module sense).

 \smallskip\noindent {\it Example.} The sequence $(*)$ splits if $R$ is a normal $\Q$-algebra (take for $\lambda$  the trace divided out by  the degree). Hence Question 1 has an easy positive answer in this case, but Question 2 is much more difficult.

 \smallskip\noindent  {\it Counter-example.}  $(*)$ does not split for $R = \Q[x,y]/(xy)$ 
  and its normalization $S$. Besides,  $(*)$  does not always split for 
 normal $\mathbb F_p$-algebras. Question 1 and therefore Question 2 have negative answers in these cases.

\subsection{} 
 {M. Hochster's {\it direct summand conjecture} \cite{H1} (1969) concerns Question $1$, namely:
  
\smallskip DSC:  {\it $\;\;\;\;\;\;\;$ The sequence $(*)$ splits if $R$ is \emph{regular}.}

\medskip For instance, this may be the situation of a Noether normalization of an affine ring $S$ ($R$ being a polynomial ring), or a Cohen presentation of a complete local domain $S$ ($R$ being a ring of formal power series).

\medskip  The state of the art of DSC before 2016 could be summarized as follows:

$a)$  Hochster gave (short) proofs when $R$ contains a field,

$b)$ he reduced the problem to the unramified complete local case with perfect residue field $k$ of char. $p$:
$$R \cong  W(k)[[x_2, \cdots, x_d]], \;\; (x_1=p),$$
$\;\;\;\, c)$  R. Heitmann \cite{He} (2002) gave a proof in dimension $d \leq 3$\footnote{in dimension $\leq 2$, this is easy: up to replacing $S$ by its normalization, it can be assumed to be a reflexive $R$-module; it is then faithfully flat over the regular ring $R$ of dimension $\leq 2$, so that $S/R$ is finite flat, hence projective, which ensures that $(*)$ splits. The first difficult case is dimension $3$. Heitmann's proof already has a flavor of ``almost algebra", which will play a fundamental role in the general case.}.

 \medskip Here is one of Hochster's proofs for a regular local domain $R$ of char. $p$ with perfect residue field $k = R/\mathfrak m$. Let $F$ denote the Frobenius endomorphism of $R$ or its finite extension $S$, and let $n$ be a positive integer. Then $F^n \mathfrak m $ is the maximal ideal of $F^n R$ and is contained in $\mathfrak m^{p^n}$, so that $\cap \, (R.F^n \mathfrak m) = 0$ by Krull's theorem. On the other hand, since $k$ is perfect, $R$ is a finite $F^nR$-module, and since $R$ is a regular local ring, this is a free module by Kunz's theorem.  Let $\lambda$ be any nonzero $R$-linear form on $S$; rescaling, we may arrange that $\lambda(1)\neq 0$, hence $\lambda(1)\notin R.F^n \mathfrak m $ for some $n$, so that $\lambda(1)$ is part of a basis of the $F^nR$-module $R$ by Nakayama's lemma. Therefore, there exists a $F^nR$-linear form $\mu$ on $R$ such that $\mu\lambda(1)= 1$ and in particular, $\mu\lambda_{\mid F^n S}$ provides a splitting of $F^n R \subset F^n S$. By transport of structure via $F^n$, $(*)$ also splits.

\subsection{}\label{aic} As for Question 2, Hochster and C. Huneke have given a positive answer in the geometric case, \ie in presence of a base field. This is much more difficult than the above short argument, and it turns out that (unlike Question 1) the char. $0$ case is even more difficult (and is deduced from) the char. $p$ case.  They proved that if $R$ is a complete regular local domain of char. $p$ and $S$ is a finite extension domain, then the absolute integral closure $R^+$ of $R$ (\ie the integral closure in an algebraic closure of the fraction field), viewed as an $S$-algebra, is faithfully flat over $R$ \cite{HH1}.

\subsection{}\label{thms}  The direct summand conjecture is now a theorem:
  \begin{thm}[A. \cite{A2} 2016]\label{T1}
Any finite extension of a regular ring splits (as a module). 
  \end{thm}
 
 In fact, what is stronger, Question 2 also admits a positive answer for regular rings:
 
   \begin{thm}\cite{A2}\label{T2}
For any finite extension $S$ of a regular ring $R$, there is an $S$-algebra $T$ which is faithfully flat over $R$.
  \end{thm}

After Hochster's work, it suffices to deal, in both theorems, with the mixed characteristic case, and more specifically with  
  $$R= W(k)[[x_2, \cdots, x_d]]$$ where $W(k)$ stands for the Witt ring of a perfect field $k$ of char. $p$.  In that case, and if $T$ is $p$-torsion free and $p$-adically complete, {\it faithful flatness over $R$ can be checked by faithful flatness of $T/p$ over $R/p$}. 

\smallskip The simple argument sketched above for DSC in char. $p$ does not extend to the mixed characteristic case: $R$ has a natural Frobenius endomorphism, but it does not extend to the finite extension $S$ in general. This suggests to introduce a mixed characteristic analog of the perfect closure, by introducing $p^{th}$-power roots of the system of parameters $x_1 =p, x_2, \cdots, x_d $, which brings us into the perfectoid world.

\section{Perfectoid notions.}
 
\subsection{}\label{almost} Let us begin with the notion of perfectoid field, whose origin dates back to J. Tate's studies in Galois cohomology. Let $K$ be a complete non-archimedean field, 
 $K^{o}$ its valuation ring, and $K^{oo}$ its valuation ideal.

Assume that the valuation is {\it not discrete} (equivalently: $K^{oo}=  (K^{oo})^2\,)$, and that 
the residue field $k$ is of char. $p>0$.

   \medskip  
 \begin{prop}[Gabber-Ramero \cite{GR}] The following are equivalent:
 
 $i)$ 
 $K^{o}/p \stackrel{ x\mapsto x^p}{\to} K^{o}/p$ is surjective, 

 $ii)$ for each finite separable $L/K$, $L^{o}$ is {almost} etale over $K^{o}$, \ie $\Omega_{L^{o}/K^{o}}$ is killed by $K^{oo}$.\end{prop}

Such a field $K$ is called {\it perfectoid} (if one refers to $i)$) or {\it deeply ramified} (if one refers to $ii)$). The residue field $k$ is then a perfect field.

 {\it Example.} Let $K_0= W(k)[\frac{1}{p}]$ and $K$ be the completion of $K_0[p^{\frac{1}{p^\infty}}]$.  This is the basic perfectoid field in the sequel. 
 Condition $i)$ is easy to checked directly, while condition $ii)$ tells us that $K$ absorbs almost all the ramification of the completion $\hat{\bar K}_0$ of the algebraic closure of $K_0$ (this was used by Tate to compute the Galois cohomology of $\hat{\bar K}_0$).
 
\smallskip Here and in the proposition, {``almost"} is used in the sense of Almost Algebra: given a commutative ring $\mathfrak V$ and an idempotent ideal $\mathfrak m$, one ``neglects" all $\mathfrak V$-module killed by $\mathfrak m$. 
 Almost algebra, introduced by Faltings \cite{F} and developed by O. Gabber and L. Ramero \cite{GR}, goes much beyond categorical localization, and studies (non-categorical) notions such as almost finite, almost flat, almost etale.
 
\medskip We shall write {``$p^{\frac{1}{p^\infty}}$-almost"} to specify that the set-up is $(\mathfrak V,\mathfrak m) = (K^{o}, K^{oo}= p^{\frac{1}{p^\infty}}K^{o})$: ``$p^{\frac{1}{p^\infty}}$-almost zero" means ``killed by all fractional powers $p^{\frac{1}{p^i}}$".

\subsection{} The generalization from fields to algebras, and further to spaces, was initiated by G. Faltings, and fully developed by P. Scholze (and also, to some extent and independently, by K. Kedlaya and R. Liu). 
 
 Let $A$ be a Banach $K$-algebra, and 
  $A^{o}$ be its sub-$K^{o}$-algebra of power-bounded elements.

 \begin{defn}[Scholze]  {\it $A$ is  {\emph{perfectoid}} if
 $A^{o}/p^{1/p} \stackrel{ x\mapsto x^p}{\to} A^{o}/p$ is an isomorphism.} \end{defn}

    The norm of such an algebra $A$ is equivalent to the spectral norm, and $A^{o}$ is the unit ball for the latter.
   
   {\it Example.} $A^{o}  = R_\infty := \widehat{{\cup} W(k)[p^{\frac{1}{p^i}}][[x_2^{\frac{1}{p^i}}, \cdots, x_d^{\frac{1}{p^i}}]]},\; A= R_\infty[\frac{1}{p}]$.

    \medskip A fundamental result of perfectoid theory is the so-called ``almost purity theorem":  
     
\begin{thm}[Faltings; Scholze \cite{S1}, Kedlaya-Liu\cite{KL}] Let $A$ be a perfectoid algebra over a perfectoid field, and  $B$ be a finite etale $A$-algebra. 
 
   Then
   $B$ is perfectoid, and $B^{o}$ is an {  $p^{\frac{1}{p^\infty}}$-almost finite etale} $A^{o}$ algebra.  \end{thm}

\section{Perfectoid Abhyankar lemma}
 
\subsection{} Let us come back to the situation of DSC, and keep the above notation. $S\otimes_R R_\infty[\frac{1}{p}]$ may not be etale over $R_\infty[\frac{1}{p}]$, hence one cannot apply Almost Purity: a non-etale finite extension of a perfectoid algebra need not be perfectoid.   
  
\smallskip To remedy this, we take inspiration from {\it Abhyankar's classical lemma}, which tells that under appropriate assumptions, one can achieve etaleness by adjoining roots of the discriminant. 

\smallskip  We follow this strategy.  Let 
 $g\in R = W(k)[[x_2, \cdots, x_d]]$ be a discriminant of $S[\frac{1}{p}]/R[\frac{1}{p}]$. The first step is to note that adjoining $p^{th}$-power roots of $g$, in the (non-naive) sense of considering $R_\infty[\frac{1}{p}]\langle g^{\frac{1}{p^\infty}}\rangle^{o}$, is ``harmless", namely:

   \begin{thm}[A. \cite{A2}] Let $ A$ be a perfectoid $\sK$-algebra, and let $g\in  A^{o}$ be a non-zero divisor. Then for any $n$, $ A\langle g^{\e}\rangle^{o}/p^m$ is $p^{\e}$-almost faithfully flat over $ A^{o}/p^m$.
    \end{thm}
    
    In particular,  $R_{\infty, g} := R_\infty[\frac{1}{p}]\langle g^{\frac{1}{p^\infty}}\rangle^{o}$ has the property that $R_{\infty, g} /p$ is  $p^{\frac{1}{p^\infty}}$-almost faithfully flat over $R_\infty/p$.
 Note that $R_{\infty,g} $ is in general much bigger (even before completion) than the ring $R_\infty[ g^{\frac{1}{p^\infty}}]$: it contains elements such as $p^{-\frac{1}{p}}(g^{\sigma^{-1}} - g^{\frac{1}{p}})$, where $\sigma: W(k)[[x_2^{\frac{1}{p}}, \cdots, x_d^{\frac{1}{p}}]]\to W(k)[[x_2, \cdots, x_d]]$ is the natural $W(k)$-semilinear Frobenius isomorphism.

  \smallskip The proof of the theorem uses a deformation argument: one spreads out the perfectoid space attached to $A^{o}_\infty $ into a perfectoid space $Y$ by adding one variable $x$. 
 Let 
 $Y^{<\epsilon}$ be the  (perfectoid) $\epsilon$-tubular neighborhood of the locus $x=g$. Then:

 \smallskip  $\bullet$   $ A\langle g^{\e}\rangle^{o}$ is $p^{\frac{1}{p^\infty}}$-almost equal to $\; \widehat{\rm colim}_\epsilon\, \;\mathcal O^+(Y^{<\epsilon})$. 
 
  \smallskip  $\bullet$ Using  
 Scholze's (almost) description of $\,\mathcal O^+(Y^{<\epsilon})/p^\epsilon\,$ in terms of "Puiseux-like" series in the variable $x-g$ with coefficients in $\,O^+(Y)/p^\epsilon$, one shows that $\,\mathcal O^+(Y^{<\epsilon})/p^\epsilon\,$ is almost faithfully flat over  $\mathcal O^+(Y)/p^\epsilon$.
 
   For details, we refer to \cite{A2} (and/or to the short account \cite{A4})\footnote{this theorem and its avatars also play an important role in the latest development of perfectoid theory: prism and prismatic cohomology (Bhatt, Scholze).}.

 \subsection{} The next step represents at the same time a perfectoid version of Abhyankar's lemma, and a ramified version of the Almost Purity theorem.

\begin{thm}[A. \cite{A1}]  Let $A$ be a perfectoid $K$-algebra, containing a sequence of $p^{th}$-power roots of a non-zero divisor $g \in A^{o}\,$  ({\it Ex:} $A= R_{\infty, g}[\frac{1}{p}]$).
  Let $B'$ be a finite etale $A[\frac{1}{g}]$-algebra, and  
   $B^{o}$ be the integral closure of $ A^{o}$ in $B'\,$ (hence $B^{o}[\frac{1}{pg}]= B'$).  
 
 Then for any $n$, $B^{o}/ p^n$ is an {  $(pg)^{\frac{1}{p^\infty}}$-almost finite etale} $A^{o}/p^n$ algebra.
\end{thm}
 
The proof also uses a deformation argument.   
Let $X$ be the perfectoid space attached to $A^o= \mathcal O^+(X)$, and
  $X_{ >\epsilon}$  be the (perfectoid) complement of $\epsilon$-tubular neighborhood of the discriminant locus $  \; g=0$. Then:
 
 \smallskip  $\bullet$  $\mathcal O^+(X)$ is $(pg)^{\frac{1}{p^\infty}}$-almost equal to $\;   {\rm lim}_\epsilon \;\mathcal O^+(X_{  >\epsilon})$ (by Scholze's {perfectoid Riemann extension theorem}). 
 
 \smallskip  $\bullet$ By almost purity over $X_{  >\epsilon},\;(B'\mathcal O(X_{  >\epsilon}))^+$ is almost finite etale over $\mathcal O^+(X_{ >\epsilon})$. 
 
 The bulk of the work is the passage to the limit $\epsilon \to 0$, for which we refer to \cite{A1} (and/or again to the short account \cite{A4}).

\subsection{} Let us go back to DSC, in the mixed characteristic setting:

\smallskip $R= W(k)[[x_2, \cdots, x_d]], \;S\,$ finite extension, etale outside $pg = 0$.

\medskip\noindent Applying the perfectoid Abhyankar lemma with $A= R_{\infty, g}[\frac{1}{p}], \, B'= S\otimes_R A[\frac{1}{g}]$, one gets an  $S$-algebra $B^{o}$ sitting on top of a tower 
$$R \stackrel{\alpha}{\to} R_\infty \stackrel{\beta}{\to} R_{\infty, g} \stackrel{\gamma}{\to}  B^{o}$$
where $\alpha$ is faithfully flat, 
 $\beta$ is $p^{\frac{1}{p^\infty}}$-almost faithfully flat mod. $p$, 
 $\gamma$ is $(pg)^{\frac{1}{p^\infty}}$-almost faithfully flat mod. $p$. It follows that $B^{o}/p$ is almost isomorphic to a faithfully flat $R/p$-algebra.
 {\it Thus $B^{o}$ is ``almost" our wanted $\,T\,$} (\cf \ref{thms}).

\medskip {\it  
 How to get rid of ``almost"?} For this, we will use some ``Cohen-Macaulay notions".

\section{Big Cohen-Macaulay algebras}

\subsection{} Let $S$ be a Noetherian local ring, and $T$ be a (possibly big, \ie non-Noetherian) extension.

 \begin{defn}{$T$ is a {\it $\,$(big) Cohen-Macaulay $S$-algebra} if any sequence of parameters  $x_1, \ldots, x_d$ of $S$ becomes  regular in $T\,$   (\ie  $x_1$ is non-zero-divisor in $T$,  $x_2$ is non-zero divisor in $T/x_1$ ..., and $T \neq (x_1, \ldots, x_d)T$). 
} \end{defn}
 
Hochster conjectured that {\it such a $T$ always exists}. 

\smallskip 
{\it Example.} If $S$ a complete local domain of char. $p$, its absolute integral closure $S^+$ is a big Cohen-Macaulay $S$-algebra (Hochster-Huneke, \cf \ref{aic}).
 
 \smallskip The relevance of Cohen-Macaulay algebras to our purpose comes from the following (classical) lemma:

 \begin{lemma}{Assume $S$ is a complete local domain, hence a finite extension of complete regular local ring $R\,$. 
 
 Then  $T$ is a (big) Cohen-Macaulay $S$-algebra if and only if $T$ is $R$-faithfully flat.   
}\end{lemma}

Flatness of Cohen-Macaulay algebras $T$ over a regular local ring $R$ can be seen as follows. Any finite $R$-module $M$ has projective dimension $\leq d=\dim R$, hence  ${\rm{Tor}}^{R}_i(M, T)=0$ for all $i>d$. One argues by descending induction on $i$. By d\'evissage, one may assume that $M= A/\frak p$ for a prime $\frak p$. If $(x_1, \ldots, x_i)$  is a maximal regular sequence contained in $\frak p$, $M$ embeds into $N := R/ R(x_1, \ldots, x_i)$. But ${\rm{Tor}}^{R}_i(N, T)=0$ if $T$ is a Cohen-Macaulay $R$-algebra. The exact sequence ${\rm{Tor}}^{R}_{i+1}(N/M, T)\to {\rm{Tor}}^{R}_i(M, T) \to {\rm{Tor}}^{R}_i(N, T)$ is thus the $0$-sequence, hence $T$ is a flat $R$-module.

\subsection{}  The problem to pass from an {\it almost} Cohen-Macaulay $S$-algebra $T_-$ (such as $B^{o}$ at the end of the previous section) to a {\it genuine} Cohen-Macaulay $S$-algebra $T$ can be settled by Hochster's (classical) technique of algebra modifications \cite{H3}\cite{HH2}. The idea is to force each relation $x_{i+1}t_{i+1} = \sum_1^i x_j t_j$ with $t_j\in T_-, x_i\in S$, to come from a relation $t_{i+1} = \sum_1^i x'_j t'_j$ by introducing step by step new variables: a partial modification of degree $n$ of an $S$-module $M$ is a homomorphism $M\to M'$ where, given a relation $x_{i+1}m_{i+1} = \sum_1^i x_j m_j$ with coefficients $m_j$ in $M$, $M': = M[T_1, \ldots, T_i]_{\leq n}/ (m_{i+1} - \sum_1^i x_j T_j)\cdot M[T_1, \ldots, T_i]_{\leq n -1} . $ Starting from $T_-$ and taking the colimit of partial modifications for all such relations yields an $S$-algebra $T$ for which $x_{i+1}$ is not a zero divisor mod. $x_1, \ldots, x_i$. Condition $T \neq (x_1, \ldots, x_d)T$ comes from the assumption that $T_-$ is almost Cohen-Macaulay.

 \smallskip This provides the last step in the proof of Hochster's conjecture about the existence of big Cohen-Macaulay algebras. Using the above lemma, this implies Theorem \ref{T2} and the Direct Summand Conjecture \ref{T1}  (positive answer to Question 1 and 2 for regular rings).

\subsection{} This is the first item of the following

\begin{thm}[A. \cite{A2}, \cite{A3}]\label{main}   
\begin{enumerate}\item  
Any Noetherian local ring $S$ has a big Cohen-Macaulay algebra $T$.

\item  
For any local morphism $S\to S'$ of Noetherian complete local domains, there is a morphism of respective big Cohen-Macaulay algebras $T\to T'$.
\end{enumerate}
\end{thm} 

The second item is the so-called ``weak functoriality of big Cohen-Macaulay algebras", also conjectured by Hochster. Its poof is more subtle and requires a delicate consideration of integral perfectoid algebras (the most difficult case is when $S$ is of mixed characteristic and $S'$ of char. $p$).

\smallskip In mixed characteristic, one has a more precise result regarding item $(1)$ (which is used in the proof of item $(2)$):

\begin{thm}[Shimomoto \cite{Sh2}, A. \cite{A3}]\label{main'} Any Noetherian complete local domain $S$ of char. $(0,p)$ admits an (integral) {perfectoid}\footnote{an integral perfectoid $K^{o}$-algebra is a $p$-adically complete $p$-torsionfree $K^{o}$-algebra $T$ such that the Frobenius map $T/p^{\frac{1}{p}}\to T/p$ is an isomorphism.} Cohen-Macaulay algebra $T$.
 \end{thm}
  
K. Shimomoto's idea is to use the tilting equivalence between perfectoid algebras over $K$ and certain perfect algebras over a perfect field $K^\flat\cong \widehat{k((t^{\frac{1}{p^\infty}}}))$ of char. $p$ (Scholze): one can apply Hochster modifications in char. $p$, after tilting, and then untilt.  An alternative and more recent argument, due to Gabber \cite{G}, uses ultraproducts instead of tilting and Hochster modifications.

\subsection{Kunz' theorem in mixed characteristic} In the same vein but independently from the above results, B. Bhatt, S. Iyengar and L. Ma extended Kunz's flatness criterion in char. $p\,$ to mixed characteristic. 
 Let us recall that for a Noetherian ring $R$ of char. $p$, Kunz' theorem asserts that {\it $R$ is regular if and only if $R\stackrel{x\mapsto x^p}{\to} R$ is flat} (which amounts to saying that {\it there exists a {perfect}, faithfully flat $R$-algebra}\footnote{such an algebra is a (big) Cohen-Macaulay $R$-algebra.}).

\medskip Let $R$ now be any Noetherian $p$-adically complete ring.  

\begin{thm}[Bhatt-Iyengar-Ma \cite{BIM} 2018] 
  $R$ is regular if and only if there exists an integral {perfectoid}\footnote{in some generalized sense, which does not require a perfectoid base field.}, faithfully flat, $R$-algebra. \end{thm}

\section{Applications to the homological conjectures} 

\subsection{The homological turn} We have already mentioned the 
homological turn in commutative algebra in the 60's: from the 
  study of Noetherian rings and their ideals (Krull, Zariski...) to the study of the
  homological properties of their modules (Auslander, Buchsbaum, Serre...).
 
  \smallskip {\it Example.} {\it A local ring $R$ is regular if and only if every finite $R$-module has a finite free resolution. }
 
 \medskip C. Peskine and L. Szpiro introduced reduction techniques to char. $p$, and the crucial observation that in char. $p$, Frobenius preserves 
 finite free resolutions (which may be seen as an extension of Kunz' theorem).  
 
 \smallskip {\it Example.} {\it A local ring $S$ is Cohen-Macaulay if and only if there is an $S$-module of finite length with a finite free resolution \cite{PS}, \cite{Ro}. }  
 
 \smallskip We refer to \cite{BH}\cite{H4} for accounts of the homological conjectures (before the introduction of perfectoid techniques). The relevance of Cohen-Macaulay algebras typically comes from the following observation: let  $$(F_\ast)\;\; \;\; 0 \to F_s\to F_{s-1} \to  \ldots \to F_i \stackrel{\phi_i}{\to} F_{i-1} \to \ldots  F_1 \to F_0$$ be a complex of finite free $S$-modules, let $r_i$ be the rank of ${\rm{Im}}\, \phi_i$, and let $I_{r_i}(\phi_i)$ be the Fitting ideal of $S$ generated the minors of size $r_i$ of $\phi_i$. 
 
If $(F_\ast)$ is acyclic, then ${\rm{codim}}\,I_{r_i}(\phi_i)\geq i$ 
(Buchsbaum-Eisenbud). The converse is not true in general if $S$ is not Cohen-Macaulay; however if ${\rm{codim}}\,I_{r_i}(\phi_i)\geq i$, 
 then $(F_\ast)\otimes_S T $ is acyclic for any Cohen-Macaulay $S$-algebra $T $ (\cf \cite[9.1.8]{BH}).

 \subsection{Pure subrings of regular rings} The following consequence of weak functoriality of Cohen-Macaulay algebras generalizes and extends in mixed characteristic the well-known theorem of Hochster-Roberts about Cohen-Macaulayness of rings of invariants of polynomials under reductive groups in char. $0$. 
 
  \begin{thm}[Heitmann-Ma \cite{HeMa1}; A. \cite{A3}] Any pure subring of a regular ring (\eg any subring which is a direct summand as a module) is Cohen-Macaulay.\end{thm}

 Let us recall the (classical) derivation of the corollary from \ref{main} $(2)$. We denote the embedding by $S\inj R'$, where $R'$ is regular. We may also assume that $S$ is a finite extension of a regular ring $R$, and it suffices to show that $S$ is flat over $R$.  
 Let $T_S\to T_{R'}$ be a compatible map of Cohen-Macaulay algebras.  For any $R$-module $M$, one has a commutative square    
\[   \begin{CD} {\rm{Tor}}_i^{R}(M, S)  @>a>> {\rm{Tor}}_i^{R}(M, R')     \\       @V VV   @VbVV     \\\     {\rm{Tor}}_i^{R}(M, T_S) @>  >>  {\rm{Tor}}_i^{R}(M, T_{R'})   . \end{CD}\] in which $a, b  $ are injective since $S\inj R'\inj T_{R'}$ are pure morphisms (see \eg \cite[5.1.1]{A2}), and $  {\rm{Tor}}_i^{R}(M, T_S) =0$ for $i>0$ since $T_S$ is Cohen-Macaulay over $S$, hence flat over $R$; therefore ${\rm{Tor}}_i^{R}(M, S) =0$.

  In fact, Heitmann and Ma proved just the right amount of weak functoriality needed to get the corollary, but they also obtained more:  $S$ has {\it pseudo-rational singularities}. 

   \subsection{The syzygy theorem} Let us turn to the traditional problem of estimating Betti numbers and ranks of syzygies of noetherian modules.  

Let $S$ be a Noetherian local ring, and let $M$ be a finite $S$-module which admits a finite free resolution.

\begin{thm} 
   Let $0\to S^{b_s}\to S^{b_{s-1}}\to \ldots \to S^{b_0} \to M\to 0$ be a minimal free resolution of $M$. Then the $i$-th syzygy has rank $\geq i$ for any $i\leq s-1$. 
    
 A fortiori, $b_i\geq 2i+1$ if $ i<s-1$, and $b_{s-1}\geq s$. \end{thm} 
 
 These bounds are optimal in general; for instance, if $S$ is the localization of $k[x_1, x_2, x_3]$ at $0$ and $M=k$, the minimal free resolution is $0\to S \to S^{3} \to S^{3} \to S \to k\to 0$, so that $s=3$ and the first two syzygies have rank $1$ and $2$ respectively.

 E. Evans and P. Griffiths  \cite{EG} stated and proved this theorem when $S$ contains a field, using the existence of Cohen-Macaulay algebras in this case. Hochster \cite{H2} and T. Ogoma \cite{O} proved that the statement follows from DSC in general.
  DSC being proven (thanks to perfectoid techniques), these bounds are now {\it unconditional}.
  
  Similar bounds hold for minimal injective resolutions, $S$ being replaced by the injective hull of the residue field (\cf \cite[9.6]{BH}). 
 
\smallskip Actually, all homological conjectures which were standing for a while on Hochster's list are now solved, as consequences of theorem \ref{main}.

\section{Applications to singularities}

\subsection{} Starting from theorems \ref{main} and \ref{main'}, Ma and K. Schwede \cite{MaS1}\cite{MaS2} developed an analog of {tight closure theory} in mixed characteristic, with applications to singularities.
  In their theory, integral perfectoid Cohen-Macaulay algebras play somehow the role of resolution of singularities in char. $0$.  
  
  Let $S$ be a local domain, essentially of finite type over $\C$, and let $\pi: Y \to {\rm{Spec}}\, S$ be a resolution of singularities. By Grauert-Riemenschneider, $R^i\Gamma(Y, \omega_Y)= 0$ for $i>0$, whence $\mathbb H^j_{\frak m}(R\Gamma(Y, \mathcal O_Y))=0$ for $j< \dim S$ by local duality. Thus $R\Gamma(Y, \mathcal O_Y)$ appears as a ``derived avatar" of a Cohen-Macaulay algebra. Since in many questions about singularities, it is not so much resolution itself which matters but the object $R\Gamma(Y, \mathcal O_Y)\in D^b(S)$, the idea is to replace, in mixed characteristic or in char. $p$, this object by suitable (big) Cohen-Macaulay $S$-algebras. 
 
  \subsection{}  Let us just outline one striking application of this idea to {\it rational singularities}. 
    Recall that, by definition, $S$ (as before) ``is" a rational singularity if and only if $R\Gamma(Y, \mathcal O_Y) \cong S $. In particular, by Grauert-Riemenschneider and local duality, $S$ is Cohen-Macaulay. 
 
 It has been known for a long time (N. Hara \cite{Ha}, K. Watanabe, K. Smith...) that  $S$ is a rational singularity if and only if, after ``spreading out",  $(S$ mod. $p)\,$ is a $F$-rational singularity (\ie local cohomology is a simple Frobenius module) for $p>>0$. While $F$-rationality is a checkable property ({\rm Macaulay2}), checking it for all $p>>0$ does not give rise to an algorithm. 
 
 But remarkably, Ma and Schwede recently proved \cite{MaS2} that it suffices to check that  {\it $(S$ mod. $p)\,$ is a $F$-rational singularity for {some} $p$}. This provides an efficient algorithm: in practice, one may choose a very small prime $p$. 
 
 The proof goes through mixed characteristic. More precisely, it uses 
 a perfectoid avatar of the notion of rational singularity in mixed characteristic as link between char. $p$ and  char. $0$.   In some sense, this is an application of perfectoid theory to complex geometry!  
 
 A similar result holds for {\it log-terminal singularities} (rational singularities for which the multiplier ideal is trivial); in this case, the corresponding property in char. $p$ has a simple formulation: {\it for any $s\in S$, there is a positive integer $n$ such that $S\stackrel{\cdot s^{\frac{1}{p^n}}}{\to} S^{\frac{1}{p^n}}$ splits}.
 
\smallskip These results may be the beginning of a program toward making the connection between singularity theory and birational geometry in char. $0$ and their counterparts in char. $p\,$ both effective and direct (even though the proof that the algorithms work would be indirect and go through mixed characteristic and perfectoid techniques).

\subsection{} Another striking application concerns symbolic powers. For any Noetherian ring $S$, prime ideal $\frak p$, and positive integer $n$, the symbolic power is defined by
$$ \frak p^{(n)} := (\frak p^n S_{\frak p}) \cap \frak p  .$$
In the affine case, this is just the ring of functions which vanish at $V(\frak p)$ at order at least $n$. Obviously, $ \frak p^{(n)} \supset \frak p^n$, and if $\frak p$ is generated by a regular sequence, it is well-known that $ \frak p^{(n)}= \frak p^n$. To compare $ \frak p^{(n)}$ and $ \frak p^n$ in general is a classical problem which has far-reaching applications, \eg in complex analysis or in transcendental number theory (Waldschmidt constants).

 \begin{thm}[Ma, Schwede \cite{MaS1}] Let $S$ be regular of dimension $d$. Then for any prime $\frak p$ and any positive integer $n$, $ \frak p^{(dn)}\subset \frak p^n.$
 \end{thm}

In char. $0$, this was proven by Ein-Lazarsfeld-Smith using subadditivity of the multiplier ideal. In mixed characteristic, Ma and Schwede follow the same strategy, defining a new notion of multiplier ideal in which the complex  $R\Gamma(Y, \mathcal O_Y)$ attached to a log-resolution of $V(\frak p)$ is replaced by an integral perfectoid  Cohen-Macaulay algebra for $S_{\frak p}$.

     }\end{sloppypar}

    \end{document}